\input amstex
\documentstyle{amsppt}
\magnification1200
\def\w{\omega}
\def\ti{\times}
\def\IN{\Bbb N}
\def\IR{\Bbb R}
\def\pr{\operatorname{pr}}
\def\conv{\operatorname{conv}}
\def\e{\varepsilon}
\def\dlim{\varinjlim}
\def\bs{\backslash}
\def\M{\Cal M}
\def\K{\Cal K}
\def\spa{\operatorname{span}}
\def\a{\alpha}

\vsize9.2truein

\topmatter
\title
On topological groups containing a Fr\'echet-Urysohn fan
\endtitle
\author
Taras Banakh
\endauthor
\abstract{Suppose $G$ is a topological group containing a (closed) topological
copy of the Fr\'echet-Urysohn fan. If $G$ is a perfectly normal sequential
space (a normal $k$-space) then every closed metrizable subset in $G$ is
locally compact. Applying this result to topological groups whose underlying
topological space can be written as a direct limit of a sequence of closed 
metrizable subsets, we get that every such a group either is metrizable or 
is homeomorphic to the
product of a $k_\w$-space and a discrete space.
}\endabstract
%\address Department of Mathematics, Lviv University, Universytetska 1,
%Lviv, 290602, Ukaraine\endaddress
\endtopmatter

\document
The present investigation was stimulated by the paper \cite{Pe} of
E.~Pentsak who studied the topology of the direct limit $X^\infty=\dlim
X^n$ of the sequence
$$
X\subset X\ti X\subset X\ti X\ti X\subset\dots,
$$
where $(X,*)$ was a ``nice'' pointed space and $X^n$ was identified with
the subspace $X^n\ti\{*\}$ of $X^{n+1}$. In particular, in \cite{Pe} the
topology of the direct limit $l_2^\infty$, where $l_2$ is the separable
Hilbert space, was characterized. To characterize the space $l_2^\infty$,
it was necessary to glue together maps into $l_2^\infty$ and at this point
it turned out that the equiconnected function generated by the natural
convex structure on $l_2^\infty$ was discontinuous. The same concerned the
 addition operation on $l_2^\infty$ --- it was discontinuous.

So question arose: is $l_2^\infty$ homeomorphic to a topological group
or a convex set in a linear topological space? 

We pose this problem more generally: find simple conditions on a
topological space $X$ under which $X$ does not support certain algebraic
structure.

In order to answer this question, we will define spaces $K$, $V$, and $W$,
called test spaces, and will prove that an existence in $X$  (closed)
subspaces homeomorphic to one (or several) of the spaces $K$, $V$, $W$
forbids $X$ to carry certain algebraic structures.

Now we define two of three test spaces.

1) {\it The space $K$}. Let
$$
K=\{(0,0)\}\cup\{(\tfrac1n,\tfrac1{n\cdot m})\mid n,m\in\IN\}\subset\IR^2.
$$
The space $K$ is metrizable and not locally compact. Moreover, $K$ is a
minimal space with these properties in the sense that each metrizable
non-locally compact space contains a closed copy of $K$. In sake of
simplicity of denotations in the sequel, put $x_0=(0,0)$ and
$x_{n,m}=(\frac1n,\frac1{nm})$, $n,m\in\IN$. Thus $K=\{x_0,x_{n,m}\mid
n,m\in\IN\}$.
\vskip5pt

2) {\it The space $V$ (the Fr\'echet-Urysohn fan).} Let
$S_0=\{0\}\cup\{\frac1n\mid n\in\IN\}$ denote the convergent sequence and
let
$$V=\IN\ti S_0/\IN\ti\{0\}.
$$
Denote by $\pi_V:\IN\ti S_0\to V$ the quotient map. Let
$y_{n,m}=\pi_V(n,\frac1m)$, $n,m\in\IN$, and $y_0$ be the (unique)
non-isolated point of $V$. So $V=\{y_0,y_{n,m}\mid n,m\in\IN\}$.
Evidently, for every $n\in\IN$ the sequence $\{y_{n,m}\}_{m=1}^\infty$
converges to $y_0$. For each $k\in\IN$ let $V_k=\{y_0,y_{n,m}\mid n\le k,
m\in\IN\}$. It is easy to see that $V$ has the direct limit topology with
respect to the sequence $V_1\subset V_2\subset\dots$ (that is a set
$U\subset V$ is open if and only if the intersection $U\cap V_n$ is open
in $V_n$ for every $n\in\IN$).

A space $X$ contains a closed copy of $V$, provided $X$ can be written as
a direct limit of a sequence
$$
X_1\subset X_2\subset\dots,
$$
where each $X_n$ is a closed metrizable subset of $X$, nowhere dense in
$X_{n+1}$. In particular, the space $l_2^\infty$ contains a topological
copy of $V$. 
\vskip5pt

We call a subset $A$ of a topological group $G$ {\it multiplicative} if
for every $a,b\in A$ we have $a*b\in A$ (here $*$ stands for the group
operation on $G$). The following theorem implies that $l_2^\infty$ carries
no topological group structure.

\proclaim{Theorem 1} A normal $k$-space $X$ containing closed copies of
the test spaces $K$ and $V$ is homeomorphic to a) no closed multiplicative
subset of a topological group and b) no closed convex set in a linear
topological space.
\endproclaim

\demo{Proof} Assume the converse and let $f:K\ti V\to X$ be the map
defined for $(x,y)\in K\ti V\subset X\ti X$ by $f(x,y)=x*y$ if $X$ is a
closed multiplicative subset of a topological group with $*$ standing for
the group operation, or by $f(x,y)=\frac12x+\frac12y$ if $X$ is a closed convex
set in a linear topological space. It is easily verified that the map $f$
has the following properties:
%\roster
\item{1)} the map $(\pr_K,f):K\ti V\to K\ti X$ is a closed embedding;
\item{2)} the map $f_{y_0}:K\to X$ defined by $f_{y_0}:x\mapsto f(x,y_0)$,
$x\in K$, is a closed embedding.

Denote by $\conv(K)=\{(0,0)\}\cup\{(x,y)\mid 0<y\le x\le 1\}$ the convex
hull of $K$ in $\IR^2$ and let $h:X\to\conv(K)$ be a continuous extension
of the map $f_{y_0}^{-1}:f_{y_0}(K)\to K$, see \cite{Hu, p.63}. 
%(This extension $h$ can be
%constructed as follows. Observe that $\conv(K)$ being a convex
%$G_\delta$-subset of the plane, admits a closed embedding
%$i:\conv(K)\to\IR^\w$ and a retraction $r:\IR^\w\to\conv(K)$, $r\circ
%i=\id$. By Tietze Extension Theorem, the map $i\circ
%f_{y_0}^{-1}:f_{y_0}(K)\to\IR^\w$ can be extended to a map $g:X\to\IR^\w$.
%Then the map $h=r\circ g:X\to\conv(K)$ extends $f_{y_0}^{-1}$).

For $n,m\in\IN$ let $\e_{n,m}=\frac1{2nm(m+1)}$ and set
$$
O_{n,m}=\conv(K)\cap\big((\tfrac1n-\e_{n,m},\tfrac1n+\e_{n,m})\ti
(\tfrac1{nm}-\e_{n,m},\tfrac1{nm}+\e_{n,m})\big).
$$

One can check that $O_{n,m}$, $n,m\in\IN$, is a collection of pairwise
disjoint neighborhoods of the points $x_{n,m}$ in $\conv(K)$.
Since $lim_{m\to\infty}y_{n,m}=y_0$ and $h\circ
f(x_{n,m},y_0)=x_{n,m}=(\frac1n,\frac1{nm})$, for every $n,m\in\IN$, we
may find a number $k(n,m)\in\IN$ such that $h\circ
f(x_{n,m},y_{n,k(n,m)})\in O_{n,m}$. Moreover, these numbers can be chosen
so that $k(n,m+1)>k(n,m)$ for $n,m\in\IN$. Let
$$
Z=\{f(x_{n,m},y_{n,k(n,m)})\mid n,m\in\IN\}.
$$
Since $h\circ f(x_0,y_0)\notin O_{n,m}$ for all $n,m$, we get
$f(x_0,y_0)\notin Z$. We claim that $Z$ is closed in $X$. Since $X$ is
a $k$-space, it suffices to show that for every compact subset $C\subset
X$ the intersection $C\cap Z$ is closed in $C$. Observe that if
$f(x_{n,m},y_{n,k(n,m)})\in C$ then $h(C)\cap O_{n,m}\ne\emptyset$. Since
$h(C)\subset \conv(K)$ is compact, we get the set
$C_1=\{x_0\}\cup\{x_{n,m}\mid h(C)\cap O_{n,m}\ne\emptyset\}$ is compact
too. Because the map $(\pr_K,f):K\ti V\to K\ti X$ is a closed embedding,
the set $(\pr_K,f)^{-1}(C_1\ti C)\subset K\ti V$ is compact and its
projection $C_2$ onto $V$ is compact too. Since $V=\dlim V_n$, we get
$C_2\subset V_{n_0}$ for some $n_0$. Observe that $C\cap Z\subset f(C_1\ti
C_2)$ and thus $C\cap Z\subset \{f(x_{n,m},y_{n,k(n,m)})\mid n\le n_0,
x_{n,m}\in C_1\}$. Because of compactness of $C_1$ the latter set is
finite. Then $C\cap Z$ is finite and hence closed in $C$.

Therefore $Z\not\ni f(x_0,y_0)$ is a closed set in $X$. Using continuity
of $f$, find neighborhoods $U(x_0)\subset K$ and $V(y_0)\subset V$ of
$x_0$ and $y_0$ such that $f(U(x_0)\ti V(y_0))\cap Z=\emptyset$. Fix $n$
such that $x_{n,m}\in U(x_0)$ for every $m$. Since the sequence
$\{y_{n,m}\}_{m=1}^\infty$ converges to $y_0$ and the sequence
$\{k(n,m)\}_{m=1}^\infty$ is increasing we may find $m$ such that
$y_{n,k(n,m)}\in V(y_0)$. Then $f(U(x_0)\ti V(y_0))\cap Z\ni
f(x_{n,m},y_{n,k(n,m)})$ is not empty, a contradiction.\qed
\enddemo

In light of Theorem 1 the following Question arises.

\definition{Question} Is $l_2^\infty$ homeomorphic a) to  a multiplicative
subset of a topological group, b) to a convex set in a linear topological
space?
\enddefinition

We will give a negative answer to the question a) under the additional
assumption that the multiplicative subset contains an idempotent (the
unity of the group).
This follows from topological homogeneity of $l_2^\infty$ and the
following theorem.

\proclaim{Theorem 2} A perfectly normal sequential space $X$ containing a
closed copy of $K$ and a copy of $V$ is homeomorphic to a) no closed
convex set in a linear topological space, b) no closed
multiplicative subset of a topological group, c) no multiplicative subset
of a topological group such that the nonisolated point $y_0$ of $V\subset
X$ is an idempotent.
\endproclaim

\demo{Proof} Assume the converse and similarly as in the previous proof
define the map $f:K\ti V\to X$. Observe that the map $f$ has the following
properties:

{\parindent30pt
\item{1')} the map $(\pr_k,f):K\to V\to K\ti X$ is an embedding;
\item{2')} the map $f_{y_0}:K\to X$ is a closed embedding (in the case (c)
$f_{y_0}$ is the identity embedding because $y_0$ is the unity of the
group).
\item{}}

Let $h:X\to \conv(K)$ be a continuous extension of the map
$f_{y_0}^{-1}:f_{y_0}(K)\to K$ such that $h^{-1}(x_0)=f(x_0,y_0)$. (The
map $h$ can be constructed as follows. Using the perfect normality of $X$,
fix a map $\lambda:X\to[0,1]$ such that $\lambda^{-1}(0)=f_{y_0}(K)$. Let
$\tilde h:X\to\conv(K)$ be any extension of the map
$f_{y_0}^{-1}:f_{y_0}(K)\to K$ and define a map $h:X\to\conv(K)$ letting
$h(x)=\lambda(x)\cdot x_{1,1}+(1-\lambda(x))\tilde h(x)$ for $x\in K$.)

Similarly as in the previous proof define neighborhoods $O_{n,m}$ and the
set $Z\not\ni f(x_0,y_0)$. As in the proof of Theorem 1, to get a
contradiction, it suffices to show that the set $Z$ is closed in $X$.
Since the space $X$ is sequential, it is enough to verify that for every
convergent sequence $S\subset X$ the intersection $S\cap Z$ is closed in
$S$. We shall show that $S\cap Z$ is always finite. Assume on the
contrary, $S\cap Z$ is infinite.

If $f(x_0,y_0)$ is not a limit point of $S$ then $S\cap Z$ is finite
because the collection\break $\{h^{-1}(O_{n,m})\}_{n,m\in\IN}$ is discrete in
$X\bs f(x_0,y_0)$ and $S\bs f(x_0,y_0)$ is compact. So assume $f(x_0,y_0)$
is a limit point of $S$. Enumerate $S\cap Z=\{z_i\}_{i=1}^\infty$.
Evidently, the sequence $\{z_i\}_{i=1}^\infty$ converges to $f(x_0,y_0)$.
For every $i\in\IN$ find (unique) $n_i,m_i$ such that
$z_i=f(x_{n_i,m_i},y_{n_i,k(n_i,m_i)})$. Observe that the sequence
$\{x_{n_i,m_i}\}_{i=1}^\infty$ converges to $x_0$. Then the sequence
$\{(x_{n_i,m_i},z_i)\}_{i=1}^\infty$ converges to $(x_0,f(x_0,y_0))$ and lies
(together with its limit) in $(\pr_K,f)(K,V)$. Since $(\pr_K,f)$ is an
embedding and the projection $\pr_V:K\ti V\to V$ is continuous, we get the
set
$$
\multline
C_2=\{y_0\}\cup\{y_{n_i,k(n_i,m_i)}\mid
i\in\IN\}=\\
\pr_V\circ(\pr_K,f)^{-1}(\{(x_0,f(x_0,y_0))\}\cup\{(x_{n_i,m_i},z_i)
\mid i\in\IN\})\subset V
\endmultline
$$ is compact. Then $C_2\subset V_{n_0}$ for some
$n_0$ and thus the sequence $\{n_i\}$ is bounded, a contradiction with the
convergence of the sequence $\{x_{n_i,m_i}\}$.\qed
\enddemo

Therefore both Theorems 1 and 2 give us that $l_2^\infty$ is homeomorphic
to no topological group. And what about its powers $(l_2^\infty)^n$? Do
they admit a compatible group structure? It turns out that the answer
here is negative too. Observe that Theorems 1 or 2 are not applicable
because the powers of $l_2^\infty$ are not $k$-spaces. So we must think
out something new.

3. {\it The test space $W$.} We let $W$ be the direct limit of a sequence
$W_0\subset W_1\subset\dots$, where the spaces $W_n\subset K\ti V$ are
defined as follows. In $K\ti V$ let us consider the points:
$z_0=(x_0,y_0)$, $z_{n,m}=(x_{n,m},y_0)$, and
$z_{n,m,p,q}=(x_{n,m},y_{p,n+q})$, $n,m,p,q\in\IN$. Let
$W_0=\{z_0,z_{n,m}\mid n,m\in\IN\}$ and $W_p=W_{p-1}\cup\{z_{n,m,p,q}\mid
n,m,q\in\IN\}$ for $p\ge 1$. It is easy to see that for every $p\ge 1$
$W_p$ is a closed subspace of $K\ti V_p$ and $W_0$ is a nowhere dense
closed copy of $K$ in $W_0\cup(W_p\bs W_{p-1})$. On the union
$W=\bigcup_{p=0}^\infty W_p$ consider the topology of the direct limit
$\dlim W_p$ allowing a subset $U\subset W$  to be open if and only if $U\cap
W_p$ is open in $W_p$ for every $p$. Observe that a space $X$ contains a
closed copy of $W$, provided $X$ can be written as the direct dimit of a
sequence $X_0\subset X_1\subset\dots$, where each $X_n$ is a closed
metrizable subset of $X$, $X_n$ is nowhere dense in $X_{n+1}$, and $X_0$
is not locally compact. Since the space $l_2^\infty$ admits such a 
representation, it contains a copy of $W$.
\vskip5pt

Remark that each direct limit $X$ of a sequence of metrizable spaces
satisfies the following property:

{\parindent25pt\item{$(\M)$} for every map $f:Y\to X$ of a metrizable
space $Y$, every point $y\in Y$ has a neighborhood $U\subset Y$ such that
$f(U)$ admits a countable neighborhood base at $f(y)$.
\item {}}

Observe that a finite product of spaces with the property $(\M)$ enjoys
this property too. Since the space $l_2^\infty$ has the property $(\M)$
and contains a copy of the test space $W$, the following theorem implies
that for every $1\le n\le\w$ the power $(l_2^\infty)^n$ does not admit a
compatible group operation.

\proclaim{Theorem 3} A topological group containing a copy of the test
space $W$ can not be embedded into a countable product of spaces
satisfying the property $(\M)$.
\endproclaim

\demo{Proof} Suppose $W\subset X\subset \Pi_{n=1}^\infty X_n$, where
each $X_n$ has the property $(\M)$. Suppose $X$ is a topological group and 
denote by $*$ the group operation and by $e$ the unity of $X$. For each
$k\in\IN$ let $Y_k=\prod_{i=1}^kX_i$ and denote by $\pr_k:X\to Y_k$ the
projection onto the first $k$-coordinates.

By induction on $k$ we shall construct increasing number sequences
$\{n(k)\}_{k=1}^\infty$,\break $\{q(k,m)\}_{m=1}^\infty$, $k\in\IN$, such that
for every $k$ the sequence
$$
\{\pr_k(z_{n(k),m}^{-1}*z_{n(k),m,k,q(k,m)})\}_{m=1}^\infty
$$
converges in $Y_k$.

Let $n(0)=0$ and suppose that for $k-1$, the number $n(k-1)$ is known.
Define a map $f_k:W_0\ti W_k\to Y_k$ letting $f_k(x,y)=\pr_k(x^{-1}*y)$
for $(x,y)\in W_0\ti W_k$. Since the space $W_0\ti W_k$ is metrizable and
the space $Y_k$, being a finite product of the spaces $X_i$'s, has the
property $(\M)$, the point $(z_0,z_0)$ has a neighborhood $U_1\ti U_2
\subset W_0\ti W_k$ such that $f_k(U_1\ti U_2)$ has a countable
neighborhood base $\{O_m\}_{m=1}^\infty$ at $f_k(z_0,z_0)=\pr_k(e)$. Pick
$n(k)>n(k-1)$ so that $z_{n(k),m}\in U_1$ for every $m$. Since for every
$m$ the sequence $\{z_{n(k),m,k,q}\}_{q=1}^\infty$ converges to
$z_{n(k),m}$, we get the sequence
$\{\pr_k(z_{n(k),m}^{-1}*z_{n(k),m,k,q})\}_{q=1}^\infty$ converges to
$\pr_k(e)$. Thus, inductively,  for every $m$ we can find a number 
$q(k,m)>q(k,m-1)$ such that $\pr_k(z_{n(k),m}^{-1}*z_{n(k),m,k,q(k,m)})\in
O_m$. The iductive step is complete.

Consider the set 
$$Z=\{z_{n(k),m,k,q(k,m)}\mid k,m\in\IN\}\subset W$$
and notice that $Z$ is closed in $W$. Since $z_0\notin Z$, we may find
neighborhoods $U(z_0),U(e)\subset X$ of $z_0$ and $e$ such that
$(U(z_0)*U(e))\cap Z=\emptyset$. Let $k$ be such that $U(e)\supset
\pr_k^{-1}(O)$ for some neighborhood $O\subset Y_k$ of $\pr_k(e)$. We may
assume the number $k$ to be so great that $z_{n(k),m}\in U(z_0)$ for every
$m$. Find finally $m$ such that 
$$
\pr_k(z_{n(k),m}^{-1}*z_{n(k),m,k,q(k,m)})\in O.
$$
Then
$z_{n(k),m}^{-1}*z_{n(k),m,k,q(k,m)}\in U(e)$ and hence the intersection
$(U(z_0)*U(e))\cap Z\ni
z_{n(k),m}^{-1}*(z_{n(k),m}^{-1}*z_{n(k),m,k,q(k,m)})=z_{n(k),m,k,q(k,m)}$
is not empty, a contradiction.\qed
\enddemo

Now let us consider some applications of the obtained results.

\heading Structure of topological groups that are
$\M_\w$-spaces\endheading

Recall that a topological space $X$ is called a {\it $k_\w$-space} if $X$
contains a countable collection $\K$ of compact subsets of $X$ such that a
subset $U$ of $X$ is open in $X$ if and only if the intersection $U\cap K$
is closed in $K$ for every $K\in\K$ (equivalently, $X$ is a $k_\w$-space,
provided $X$ is the direct limit of a sequence of its compact subsets).

We define a  topological space $X$ to be an {\it $\M_\w$-space} if $X$
contains a countable collection $\M$ of closed metrizable subsets of $X$ 
such that a
subset $U$ of $X$ is open in $X$ if and only if the intersection $U\cap M$
is closed in $M$ for every $M\in\M$ (equivalently, $X$ is an $\M_\w$-space,
if $X$ is the direct limit of a sequence of its closed metrizable 
subsets).

It turns out that an existence of a compatible group structure imposes very
strict restrictions on the topology of $\M_\w$-spaces.

\proclaim{Theorem 4} Suppose a topological group $X$ is an $\M_\w$-space.
If $X$ is not metrizable, then
\roster
\item $X$ contains a closed copy of the Fr\'echet-Urysohn fan;
\item each closed metrizable subset of $X$ is locally compact;
\item $X$ contains an open subgroup $H$ that is a $k_\w$-space;
\item $X$ is homeomorphic to a product of a $k_\w$-space and a
discrete space;
\item $X$ is homeomorphic to an open subset of a $k_\w$-space.
\endroster
\endproclaim

\demo{Proof} Suppose $X$ is not metrizable and let $e$ denote the unity of the 
group $X$. Write $X=\dlim X_n$ be the direct
limit of a sequence $\{e\}=X_0\subset X_1\subset X_2\subset\dots$ consisting
of closed metrizable subsets of $X$.
 To prove 1) we will
show that for every $n$ there is $m$ such that $e$ is a limit point of the
set $X_m\bs X_n$ in $X_m$. Fix $n$ and a decreasing neighborhood base
$\{U_i\}_{i=1}^\infty$ of $e$ in $X_n$. Since $X$ is not metrizable, each
$U_i$ is not open in $X$, and thus $U_i$ is not open in some $X_{m(i)}$.
Consequently, there is a sequence $\{y_{ij}\}_{j=1}^\infty\subset
X_{m(i)}\bs X_n$ convergent to a point $x_i\in U_i$. Let $k(i)=\min\{k\in
\IN\mid \forall j_0\in\IN\;\exists j\ge j_0$ such that $y_{ij}\in X_k\}$.
Passing to a subsequence, if necessary, we may assume that
$\{y_{ij}\}_{j=1}^\infty\subset X_{k(i)}\bs X_{k(i)-1}$.

If $m=\sup\{k(i)\mid i\in\IN\}<\infty$ then all the points $y_{ij}$,
$i,j\in\IN$ lie in the set $X_m\bs X_n$. Since $X_m$ is metrizable and
the sequence $\{x_i\}$ tends to $e$, we may choose a subsequence
$\{z_j\}_{j=1}^\infty\subset\{y_{ij}\mid i,j\in\IN\}$ convergent to $e$.
Thus $e$ is a limit point of the set $X_m\bs X_n$ and we are done.

Now suppose $\sup\{k(i)\mid i\in\IN\}=\infty$. Using the continuity of the
multiplication $*$, find $p\in\IN$ such that $U_p* U_p\subset X_k$ for 
some $k$. Let $i$ be such that $k(i)>k$ and $i\ge p$. Obviously, the
sequence $\{x_i^{-1}*y_{ij}\}_{j=1}^\infty$ converges to $e$. We claim
that there exists $j_0\in\IN$ such that $x_i^{-1}*y_{ij}\notin X_n$ for
all $j\ge j_0$. Assuming the converse we would find $j$ such that
$x_i^{-1}*y_{ij}\in U_p\subset X_n$. Then $y_{ij}\in x_i*U_p\subset U_p*U_p
\subset X_k$, a contradiction with $k(i)>k$ and $y_{ij}\in X_{k(i)}\bs
X_{k(i)-1}$. Thus we have proven that $\{x_i^{-1}*y_{ij}\}_{j\ge
j_0}\subset X\bs X_n$ for some $j_0$. Since this sequence converges to $e$,
it is contained in some $X_m$.

Now we are ready to construct a closed copy of $V$ in $X$. Applying the
statement proved above, we may construct inductively an increasing number
sequence $\{m(i)\}_{i=1}^\infty$ and sequences $\{y_{ij}\}_{j=1}^\infty$,
$i\in\IN$, such that 
$$
\lim_{j\to\infty}y_{ij}=e,\; y_{ij}\in X_{m(i)}\bs X_{m(i-1)}, \; j\in\IN.
$$
Evidently, the set $\{e\}\cup\{y_{ij}\mid i,j\in\IN\}$ is a closed copy of
$V$ in $X$. Hence a) is proven.

To prove 2), notice that $X$, being a direct limit of a sequence of
metrizable spaces, is a perfectly normal sequential space. Since $X$
contains a copy of the test space $V$, by Theorems 1 and 2, $X$ contains
no closed copy of the test space $K$. Because every metrizable non-locally
compact space contains a closed copy of $K$, every closed metrizable
subset in $X$ must be locally compact.

To prove 3), let us firstly construct an open separable subset $U\subset
X$. By 2), each $X_n$ is locally compact. Thus, we may choose inductively
open neighborhoods $U_n$ of $e$ in $X_n$ so that the closure $\bar U_n$ is
compact and $\bar U_n\subset U_{n+1}$ for every $n\in\IN$. Since each
$\bar U_n$ is a (separable) metric compactum, the union
$U=\bigcup_{n\in\IN}U_n$ is an open separable subset of $X$. Then its span
$H=\spa(U)$ is an open separable subgroup of $X$. Using the separability
of $H$, one can easily show that for every $n$ the locally compact space
$H_n=H\cap X_n$ is separable, and thus $H_n$ is a $k_\w$-space. Then
$H=\dlim H_n$ is a $k_\w$-space too. Since the subgroup $H$ is open in $X$ the
decomposition of $X$ onto right residue classes of $H$ just provides us with a
homeomorphism of $X$ onto the product $H\ti D$ for some discrete space
$D$. Let $\a D$ be the one-point compactification. Evidently, $H\ti D$ is an 
open subset of the $k_\w$-space $H\ti\a D$, thus (5) follows.
\qed
\enddemo

\enddocument

\Refs\widestnumber\key{PP}

\ref\key{Pe}\by E.~Pentsak\paper On manifolds modeled on direct limits of 
 $\Cal C$-universal ANR's\jour Matematychni Studii\vol 5\year1995\pages
107--116\endref
\endRefs
\bye